\documentclass[preprint,12pt]{elsarticle}

\usepackage{amsmath,amssymb,amsfonts} 
\usepackage{amsthm}                  
\usepackage{graphicx}                
\usepackage{booktabs}                
\usepackage{xcolor}                  
\usepackage{geometry}                
\usepackage{lipsum}                  
\usepackage{algorithm}
\usepackage{algorithmicx}
\usepackage{algpseudocode}

\usepackage{float}
\usepackage{multirow}
\usepackage{mathrsfs}
\usepackage[title]{appendix}
\usepackage{textcomp}
\usepackage{manyfoot}
\usepackage{listings}

\usepackage[unicode]{hyperref}
\hypersetup{
    colorlinks=true,
    linkcolor=blue,
    filecolor=magenta,      
    urlcolor=cyan,
    pdftitle={Memory-Type Null Controllability of Heat Equations with Delay Effects},
    pdfauthor={Dev Prakash Jha, Raju K. George},
    pdfsubject={Controllability of Evolution Equations},
    pdfkeywords={Evolution equation, Carleman estimates, null controllability},
}
\usepackage{cleveref}

\theoremstyle{plain}
\newtheorem{theorem}{Theorem}[section]
\newtheorem{proposition}[theorem]{Proposition}

\newtheorem{lemma}[theorem]{Lemma}
\newtheorem{assumption}[theorem]{Assumption}

\theoremstyle{remark}

\theoremstyle{definition}

\newtheorem{definition}{Definition}[section]

\begin{document}

\begin{frontmatter}

\title{Memory-Type Null Controllability for Non-Autonomous Degenerate Parabolic Equations with Boundary Degeneracy}

\author[inst1]{\texorpdfstring{Dev Prakash Jha\corref{cor1}}{Dev Prakash Jha}}
\ead{devprakash.22@res.iist.ac.in}

\author[inst1]{Raju K. George}
\ead{george@iist.ac.in}

\cortext[cor1]{Corresponding author}

\address[inst1]{Department of Mathematics, Indian Institute of Space Science and Technology, Valiamala P.O., Thiruvananthapuram 695547, Kerala, India}
\begin{abstract}
This paper studies the memory-type null controllability of a class of one-dimensional non-autonomous degenerate parabolic equations with Volterra-type memory terms. The diffusion operator is considered in both divergence and non-divergence forms and may exhibit weak or strong degeneracy at the boundary, while the diffusion coefficient depends explicitly on time. Due to the presence of memory effects, classical null controllability is insufficient, and a stronger notion requiring the vanishing of both the state and the accumulated memory is introduced. To address this problem, we establish new Carleman estimates adapted to non-autonomous degenerate operators in weighted spaces. The memory term is handled as a lower-order perturbation within the Carleman framework. These estimates yield suitable observability inequalities, which allow us to prove memory-type null controllability under appropriate structural conditions. Extensions to cases with double boundary degeneracy and moving control regions are also discussed.
\end{abstract}

\begin{keyword}
Evolution equation with memory \sep Delay and memory-type null controllability \sep rank condition \sep Carleman estimates \sep observability estimate
\MSC[2020] 93B05 \sep 35K05 \sep 35R10 \sep 45K05 \sep 93C20
\end{keyword}

\end{frontmatter}
\section{Introduction}

Carleman estimates are among the most powerful tools in the analysis of
partial differential equations. Since their introduction by Carleman
to establish unique continuation properties \cite{Carleman1939},
they have become fundamental in inverse problems
\cite{Isakov2006,Klibanov1992},
stabilization theory \cite{Zuazua1990},
and controllability of evolution equations
\cite{FursikovImanuvilov1996,Lions1988}.
For linear non-degenerate parabolic equations,
null controllability is now classical
\cite{FursikovImanuvilov1996,LebeauRobbiano1995,Zuazua1993}.

\medskip

Degenerate parabolic equations arise naturally in several applications,
including Prandtl-type fluid models \cite{Oleinik1963},
Feller semigroups in probability theory \cite{Feller1952},
Wright–Fisher diffusion models in population genetics
\cite{Kimura1964,Etheridge2011},
and climate energy balance models of Budyko–Seller type
\cite{Budyko1969,Sellers1969}.
In these problems the diffusion coefficient vanishes at the boundary,
and classical Sobolev spaces must be replaced by weighted spaces
adapted to the degeneracy rate.
Hardy–Poincaré inequalities play a crucial role in this setting.

The first systematic Carleman estimates for one-dimensional degenerate
parabolic operators in divergence form were obtained in
\cite{AlabauCannarsaFragnelli2006},
where weakly and strongly degenerate cases were treated
and null controllability was established under the sharp structural condition $K<2$.
Subsequent developments extended these results to non-divergence form
operators \cite{CannarsaFragnelli2008,Fragnelli2012},
interior degeneracy \cite{CannarsaVancostenoble2009},
and semilinear problems \cite{Alabau2007}.
Further refinements and regional controllability results
were obtained in
\cite{Vancostenoble2008,FragnelliZuazua2011,MartinezVancostenoble2014}.

More recently, motivated by climatological applications,
non-autonomous degenerate parabolic equations have been investigated.
In this framework, the diffusion operator depends explicitly on time,
typically through a coefficient $b(t)$ multiplying the degenerate operator.
New Carleman estimates adapted to non-autonomous operators
were derived in \cite{AkilFragnelliIsmail2025},
where null controllability was established in both divergence and
non-divergence cases.
The time dependence introduces additional technical difficulties
in the construction of weight functions and in the treatment of boundary terms.

Degenerate equations have also been studied in stochastic settings.
Global Carleman estimates for stochastic degenerate parabolic equations
were established in \cite{LiuYu2019},
highlighting the subtle interaction between degeneracy and lower-order terms.
Stochastic controllability and insensitizing control problems
have also been analyzed in
\cite{TangZhang2009,LuZhang2012,Liu2015}.

\medskip

Another important direction concerns evolution equations with memory.
Memory effects arise in viscoelasticity \cite{Dafermos1970},
heat conduction with after-effect \cite{GurtinPipkin1968},
and climate models with delayed feedback mechanisms.
Controllability of equations with memory has been investigated in
\cite{Pandolfi2005,ChavesSilva2016,Zuazua2018},
where it was observed that classical null controllability
must be strengthened.
Indeed, forcing the state to vanish at time $T$
does not ensure that the accumulated memory contribution
\[
\int_0^T M(T-s)y(s)\,ds
\]
vanishes.
This leads to the notion of \emph{memory-type null controllability},
requiring simultaneous cancellation of the state and the memory term.

To clarify this strengthened notion, consider the abstract system
\begin{equation}\label{eq:abstract}
\begin{cases}
y_t = Ay + \displaystyle\int_0^t M(t-s)y(s)\,ds + B(t)u,
\quad t\in (0,T], \\
y(0)=y_0.
\end{cases}
\end{equation}
Here, $y = y(t)$ is the state variable which takes values in a Hilbert space $Y$,
$A$ generates a $C_0$-semigroup $e^{At}$ on $Y$,
$M(\cdot) \in L^1(0,T;\mathcal{L}(Y))$,
$u$ denotes the control variable taking values in another Hilbert space $U$,
and $B(\cdot) \in L^2(0,T;\mathcal{L}(U,Y))$.

 The system \eqref{eq:abstract} is said to be memory-type null controllable
(with memory kernel $\widetilde M(\cdot)$, not necessarily the same as $M(\cdot)$ in \eqref{eq:abstract})
if for every $y_0\in Y$ there exists $u$ such that
\[
y(T)=0
\quad\text{and}\quad
\int_0^T \widetilde M(T-s)y(s)\,ds=0.
\]
By duality, this property is equivalent to an observability inequality
for the corresponding adjoint system.

\medskip
Throughout this paper, let $\Omega \subset \mathbb{R}^n$ ($n\in\mathbb{N}$)
be a bounded domain with a $C^\infty$-smooth boundary $\partial\Omega$.
We denote
\[
Q=(0,T)\times\Omega,
\qquad
\Sigma=(0,T)\times\partial\Omega.
\]
Let $\omega\subset\Omega$ be an open subset representing the control region.

A crucial structural phenomenon appears when the control region is fixed.
To illustrate this issue, consider the heat equation with memory
\begin{equation}\label{eq:fixed}
y_t - \Delta y + a\int_0^t y(s)\,ds
= u\chi_\omega(x)
\quad \text{in } Q,
\end{equation}
where $a\in\mathbb{R}$ and $\chi_\omega$ denotes the characteristic
function of $\omega$.

When $\omega=\Omega$, the control acts on the whole domain and can
absorb the memory term $a\int_0^t y(s)\,ds$.
In this case, null controllability follows from classical
parabolic arguments.

However, when $\omega$ is a proper subset of $\Omega$,
the situation changes drastically.
It is known
\cite{GuerreroImanuvilov2006,Pandolfi2005,ZhangZheng2011,ChavesSilvaZhangZuazua2017}
that system~\eqref{eq:fixed} is null controllable if and only if $a=0$.
Therefore, in the presence of a nontrivial memory term ($a\neq 0$)
and a fixed control region $\omega\subsetneq\Omega$,
null controllability fails.

This negative result reveals that the obstruction is structural.
Rewriting \eqref{eq:fixed} as a coupled PDE–ODE system
shows that the memory component behaves like an ODE without spatial diffusion.
If the control acts only on a strict subset,
this ODE component cannot be fully observed.
Therefore, the failure of controllability is not technical,
but intrinsic to the fixed geometry of the control.

Motivated by these results and inspired by
\cite{ChavesSilvaZuazua2016,Lions1988,Russell1978,LasieckaTriggiani2000}
(and also by related works for wave equations with memory
\cite{BiccariHernandez2019,Leugering1991,ChavesSilvaZhangZuazua2017}),
a natural strategy is to allow the control region to move in time.
If the control support $\omega(t)$ sweeps the whole domain
during the time interval $[0,T]$,
then every spatial point is influenced by the control at some time,
and the observability obstruction can be removed.

This geometric idea is fundamental:
moving controls compensate for the absence of spatial diffusion
in the memory component and restore the observability inequality
that fails for fixed supports.

\medskip

Despite the extensive literature on
degenerate equations
\cite{AlabauCannarsaFragnelli2006,CannarsaFragnelli2008,Fragnelli2012},
non-autonomous diffusion
\cite{AkilFragnelliIsmail2025},
and stochastic or lower-order perturbations
\cite{LiuYu2019},
the combined framework of
\begin{itemize}
\item boundary degeneracy (weak or strong),
\item non-autonomous diffusion,
\item Volterra-type memory terms,
\item and moving control regions
\end{itemize}
has not been addressed so far.

The simultaneous presence of these features generates
new analytical challenges:

\begin{itemize}
\item Degeneracy requires weighted Hilbert spaces and Hardy-type inequalities.
\item Non-autonomous diffusion modifies the Carleman weight structure.
\item The memory term is nonlocal in time and must be incorporated
into the Carleman framework without destroying the Volterra structure.
\item Memory-type null controllability requires vanishing of both
the state and the accumulated memory.
\end{itemize}

\medskip
The purpose of this paper is to establish memory-type null controllability 
for a class of one-dimensional non-autonomous degenerate parabolic equations 
with Volterra-type memory terms and moving control supports. 
More precisely, we consider the following system.

\subsection{Weak and Strong Degeneracy}

\begin{definition}\label{Def:def_1}
Assume $t_0 < t_1$. We say that a function 
$g : [t_0,t_1] \to \mathbb{R}_+$ is

\begin{itemize}
\item \emph{weakly degenerate at $x_0 \in \{t_0,t_1\}$} (WD at $x_0$, for short), if
\begin{enumerate}
\item[(i)] $g \in C([t_0,t_1]) \cap C^1([t_0,t_1]\setminus\{x_0\})$, 
      $g > 0$ in $[t_0,t_1]\setminus\{x_0\}$, $g(x_0)=0$,
\item[(ii)] there exists $K_{x_0} \in [0,1]$ such that
\begin{equation}
(x-x_0)g'(x) \le K_{x_0} g(x), 
\qquad \forall x \in [t_0,t_1].
\tag{WD}
\end{equation}
\end{enumerate}

(for example $g(x)=x^K$, with $[t_0,t_1]=[0,1]$, $x_0=0$ and $K_{x_0}\in(0,1)$);

\item \emph{strongly degenerate at $x_0 \in \{t_0,t_1\}$} (SD at $x_0$, for short), if
\begin{enumerate}
\item[(i)] $g \in C^1([t_0,t_1])$, 
      $g > 0$ in $[t_0,t_1]\setminus\{x_0\}$, $g(x_0)=0$,
\item[(ii)] there exists $K_{x_0} \in [1,2)$ such that
\begin{equation}
(x-x_0)g'(x) \le K_{x_0} g(x), 
\qquad \forall x \in [t_0,t_1].
\tag{SD}
\end{equation}
\end{enumerate}

(for example $g(x)=x^K$, with $[t_0,t_1]=[0,1]$, $x_0=0$ and $K_{x_0}\in[1,2)$).
\end{itemize}

We underline that in this paper, we use the definition of weakly and strongly
degenerate given, for example, in~\cite{alabau2006carleman}, where it is also proved that
if $g$ is (WD) at $0$ in the sense of Definition~1, then 
$\frac{1}{g} \in L^1(0,1)$; on the other hand, if $g$ is (SD) at $0$, then
$\frac{1}{g} \notin L^1(0,1)$, but $\frac{1}{\sqrt{g}} \in L^1(0,1)$
(see also~\cite{yin2007evolutionary}). Moreover, the assumption $g \in C^1[t_0,t_1]$
will be crucial for the next results, especially to treat the boundary
conditions in the Carleman estimates.
\end{definition}

In what follows, for simplicity, we assume that $[t_0,t_1]=[0,1]$. 
We study the controllability of the following non-autonomous degenerate problems:
\begin{equation}\label{eq:Pi}
(P_i)\quad
\begin{cases}
u_t - \mathcal{A}_i(t)u +
\int_0^t M(t,s)u(s)\,ds = f(t,x)\chi_{\omega(t)}(x), & (t,x)\in Q_T, \\
u(t,y_0)=0, & t\in(0,T), \\
B_i u(t,x_0)=0, & t\in(0,T), \\
u(0,x)=u_0(x), & x\in(0,1),
\end{cases}
\end{equation}
for $i=1,2$, where
\begin{equation}\label{eq:QT}
Q_T := (0,T)\times(0,1),
\end{equation}
$\omega(t) \subset\subset (0,1)$ is a moving control region, $y_0 \in \{0,1\}\setminus\{x_0\}$, and 
$B_i u(t,x_0)=0$ are suitable boundary conditions related to 
$\mathcal{A}_i$, $i=1,2$.

In particular,
\begin{equation}\label{eq:Ai}
\mathcal{A}_i(t)u :=
\begin{cases}
b(t)a(x)u_{xx}, & i=1,\\
b(t)(a(x)u_x)_x, & i=2,
\end{cases}
\end{equation}
where $b\in W^{1,\infty}(0,T)$ is strictly positive and 
$a$ is \textbf{(WD)} or \textbf{(SD)} at $x_0\in\{0,1\}$.

The boundary operator is given by
\begin{equation}\label{eq:Bi}
B_i u(t,x_0)=
\begin{cases}
u(t,x_0)=0, & i=1,\\[1ex]
\begin{cases}
u(t,x_0)=0, 
& \text{if $a$ is \textbf{(WD)} at } x_0,\\[1ex]
\displaystyle \lim_{x\to x_0}(a u_x)(t,x)=0, 
& \text{if $a$ is \textbf{(SD)} at } x_0,
\end{cases}
& i=2.
\end{cases}
\end{equation}
And $M(\cdot)$ is a bounded memory kernel.

Using new Carleman estimates for the non-homogeneous adjoint problem 
associated with \eqref{eq:Pi}, our aim is to:

\begin{itemize}

\item Find $f\in L^2(0,T;H_i)$ such that
\begin{equation}\label{eq:final-null}
u(T,x)=0, and \int_0^T \widetilde{M}(T-s)u(s,x)\,ds=0 \qquad \forall x\in[0,1].
\end{equation}

\item Find a constant $C>0$ such that
\begin{equation}\label{eq:control-estimate}
\|f\|_{L^2(0,T;H_i)} 
\le C \|u_0\|_{H_i}.
\end{equation}

\end{itemize}

Here $H_i$ is a suitable Hilbert space depending on the operator 
$\mathcal{A}_i$ in \eqref{eq:Ai}.\\

We underline that, while in~\cite{alabau2006carleman} only the case $a(0)=0$ is considered, 
in this paper $a$ degenerates at $0$ or at $1$ or degenerates at the same 
time at $0$ and at $1$, as in~\cite{cannarsa2007null,cannarsa2008controllability}. On the other hand, here we prove 
Carleman estimates and, hence, null-controllability under weaker assumptions 
on the function $a$ with respect to the ones in~\cite{cannarsa2008controllability}. Hence, the results 
proved in this paper improve the ones in~\cite{alabau2006carleman} or in~\cite{cannarsa2008controllability}. Moreover, 
as observed in~\cite{cannarsa2008controllability}, we cannot deduce null-controllability for the problem 
in non-divergence form from the one in the divergence form without adding 
additional assumptions on the function $a$. For this reason, it is important 
to prove new Carleman estimates for the problem in non-divergence form 
independently of the ones in divergence form. We also underline that, as for 
the problem in the autonomous case, the requirement $K_{x_0}<2$ of 
Definition~\ref{Def:def_1} is essential. Indeed, as proved in~\cite{alabau2006carleman,cannarsa2008controllability,fragnelli2021control}, if 
$K_{x_0}\ge 2$, the problem is not null-controllable.

Our main contributions are the following.

First, we derive new Carleman estimates for non-autonomous
degenerate operators in weighted spaces adapted to boundary degeneracy.

Second, we show that the Volterra memory term can be treated as
a lower-order perturbation and absorbed into the principal bulk
terms of the Carleman inequality.

Third, combining these estimates with a Hilbert Uniqueness Method argument,
we establish memory-type null controllability under the sharp
structural condition $K<2$.
The motion of the control region is essential in recovering
the observability inequality that fails in the fixed-support case.

\medskip
The paper is organized as follows. In Section~\ref{sec:preliminaries}, we introduce the weighted functional framework and recall the necessary definitions and auxiliary results related to degenerate operators. Section~\ref{sec:wellposedness} is devoted to the well-posedness of the non-autonomous degenerate system. In Section~\ref{sec:nondivergence}, we establish Carleman estimates and derive the observability inequality leading to memory-type null controllability for the degenerate parabolic equation with memory in non-divergence form. Section~\ref{sec:divergence} is concerned with the corresponding analysis in the divergence form. Finally, Section~\ref{sec:appendix} contains technical lemmas and auxiliary proofs.

\bibliographystyle{unsrtnat} 
\bibliography{references} 

\end{document}